\documentclass[a4paper]{article}

\usepackage{amsmath,amssymb,amsthm}

\numberwithin{equation}{section}
\allowdisplaybreaks[4]

\theoremstyle{plain}

\theoremstyle{definition}
\newtheorem{definition}{Definition}[section]
\newtheorem{example}[definition]{Example}
\newtheorem{remark}[definition]{Remark}

\usepackage{graphicx}      

\usepackage{bbm}
\usepackage{mathrsfs}

\newcommand{\IZ}{\mathbbm{Z}}
\newcommand{\IR}{\mathbbm{R}}

\def\%#1{\mathcal{#1}}
\newcommand{\law}{\mathscr{L}}

\def\dtv{\mathop{d_{\mathrm{TV}}}}

\def\tv{{\mathrm{TV}}}

\newcommand{\lito}{\mathrm{o}}

\newcommand{\nsig}{\Sigma\kern-0.5em\raise0.2ex\hbox to 0pt{$\mid$}\kern0.5em}

\newcommand{\toinf}{\to\infty}

\newcommand{\eps}{\varepsilon}
\renewcommand{\phi}{\varphi}
\newcommand{\D}{\Delta}

\newcommand{\ahalf}{{\textstyle\frac{1}{2}}}

\def\a#1th{{\textstyle\frac{1}{#1}}}

\newcommand{\eq}{\eqref}
\newcommand{\IE}{\mathbbm{E}}
\newcommand{\IP}{\mathbbm{P}}

\newcommand{\Po}{\mathop{\mathrm{Po}}}

\def\be#1\ee{\begin{equation*}#1\end{equation*}}
\def\ben#1\ee{\begin{equation}#1\end{equation}}

\def\bes#1\ee{\begin{equation*}\begin{split}#1\end{split}\end{equation*}}
\def\besn#1\ee{\begin{equation}\begin{split}#1\end{split}\end{equation}}

\def\bg#1\ee{\begin{gather*}#1\end{gather*}}
\def\bgn#1\ee{\begin{gather}#1\end{gather}}

\def\bm#1\ee{\begin{multline*}#1\end{multline*}}
\def\bmn#1\ee{\begin{multline}#1\end{multline}}

\def\ba#1\ee{\begin{align*}#1\end{align*}}
\def\ban#1\ee{\begin{align}#1\end{align}}

\def\bea#1\ee{\begin{eqnarray*}#1\end{eqnarray*}}
\def\bean#1\ee{\begin{eqnarray}#1\end{eqnarray}}

\def\bklr#1{\bigl(#1\bigr)}

\def\bbbklr#1{\biggl(#1\biggr)}

\def\bkle#1{\bigl[#1\bigr]}

\def\klg#1{\{#1\}}
\def\bklg#1{\bigl\{#1\bigr\}}

\def\bbbklg#1{\biggl\{#1\biggr\}}

\def\norm#1{\Vert#1\Vert}

\def\bbbnorm#1{\biggl\Vert#1\biggr\Vert}

\def\abs#1{\vert#1\vert}
\def\babs#1{\bigl\vert#1\bigr\vert}
\def\bbabs#1{\Bigl\vert#1\Bigr\vert}

\def\mid{\vert}

\def\floor#1{{\lfloor#1\rfloor}}

\def\^#1{\ifmmode{\mathaccent"705E #1}\else{\accent94 #1}\fi}
\def\~#1{\ifmmode{\mathaccent"707E #1}\else{\accent"7E #1}\fi}
\def\*#1{#1^\ast}
\edef\-#1{\noexpand\ifmmode {\noexpand\bar{#1}} \noexpand\else \-#1\noexpand\fi}
\edef\_#1{\noexpand\ifmmode {\noexpand\underline{#1}} \noexpand\else
\underline{#1}\noexpand\fi}
\def\>#1{\vec{#1}}
\def\.#1{\dot{#1}}

\def\leq{\leqslant}

\def\geq{\geqslant}

\def\atop{\@@atop}

\newcount\minute
\newcount\hour
\newcount\hourMins
\def\now{%
\minute=\time%
\hour=\time \divide \hour by 60%
\hourMins=\hour \multiply\hourMins by 60%
\advance\minute by -\hourMins%
\zeroPadTwo{\the\hour}:\zeroPadTwo{\the\minute}%
}
\def\zeroPadTwo#1{\ifnum #1<10 0\fi#1}

\def\blfootnote{\xdef\@thefnmark{}\@footnotetext} 

\begin{document}

\title{On the Optimality of Stein Factors}
\author{Adrian R\"ollin}
\date{National University of Singapore}
 
\maketitle

\begin{abstract}The application of Stein's method for distributional approximation
often involves so called \emph{Stein factors} (also called \emph {magic
factors}) in the bound of the solutions to Stein equations. However, in some
cases these factors contain additional (undesirable) logarithmic terms. It has
been shown for many Stein factors that the known bounds are sharp and thus that
these additional logarithmic terms cannot be avoided in general. However, no
probabilistic examples have appeared in the literature that would show that
these terms in the Stein factors are not just unavoidable artefacts, but that
they are there for a good reason. In this article we close this gap by
constructing such examples. This also leads to a new interpretation of the
solutions to Stein equations.
\end{abstract}

\section{Introduction}

Stein's method for distributional approximation, introduced by \cite{Stein1972},
has been used to obtain bounds on the distance between probability measures for
a variety of distributions in different metrics. There are two main steps
involved in the implementation of the method. The first step is to set up the
so-called \emph{Stein equation}, involving a \emph{Stein operator}, and to
obtain bounds on its solutions and their derivatives or differences; this can be
done either analytically, as for example \cite{Stein1972}, or by means of the
probabilistic method introduced by \cite{Barbour1988}. In the second step one
then needs bound the expectation of a functional of the random variable under
consideration. There are various techniques to do this, such as the local
approach by \cite{Stein1972} and \cite{Chen2004a} or the exchangeable pair
coupling by \cite{Stein1986}; see \cite{Chen2009} for a unification of these and
many other approaches.

To successfully implement the method, so-called \emph{Stein factors} play an
important role. In this article we will use the term \emph{Stein factor} to
refer to the asymptotic behaviour of the bounds on the solution to the Stein
equation as some of the involved parameters tend to infinity or zero. Some of
the known Stein factors are not satisfactory, because they contain terms which
often lead to non-optimal bounds in applications. Additional work is then
necessary to circumvent this problem; see for example \cite{Brown2000}. There
are also situations where the solutions can grow exponentially fast, as has been
shown by \cite{Barbour1992c} and \cite{Barbour1998a} for some specific compound
Poisson distributions, which limits the usability of Stein's method in these
cases.

To make matters worse, for many of these Stein factors it has been shown that
they cannot be improved; see \cite{Brown1995}, \cite{Barbour1992c} and
\cite{Barbour2005}. However, these articles do not address the question whether
the problematic Stein factors express a fundamental ``flaw'' in Stein's method
or whether there are examples in which these additional terms are truly needed
if Stein's method is employed to express the distance between the involved
probability distributions in the specific metric.

The purpose of this note is to show that the latter statement is in fact true.
We will present a general method to construct corresponding probability
distributions; this construction not only explains the presence of problematic
Stein factors, but also gives new insight into Stein's method.

In the next section, we recall the general approach of Stein's method in the
context of Poisson approximation in total variation. Although in the univariate
case the Stein factors do not contain problematic terms, it will demonstrate the
basic construction of the examples. Then, in the remaining two sections, we
apply the construction to the multivariate Poisson distribution and Poisson
point processes, as in these cases the Stein factors contain a logarithmic term
which may lead to non-optimal bounds in applications.

\section{An illustrative example}

In order to explain how to construct examples which illustrate the nature of
Stein factors and also to recall the basic steps of Stein's method, we start
with the Stein-Chen method for univariate  Poisson approximation (see
\cite{Barbour1992}).

Let the total variation distance between two non-negative, integer-valued random
variables $W$ and $Z$ be defined as 
\ben                       					\label{1} 
	\dtv\bklr{\law(W),\law(Z)} 
    := \sup_{h\in\%H_{\tv}}\babs{\IE h(W) - \IE h(Z)},
\ee
where the set $\%H_{\tv}$ consists of all indicator functions on the
non-negative integers~$\IZ_+$. Assume now that $Z\sim\Po(\lambda)$. Stein's idea
is to replace the difference between the expectations on the right hand side of
\eq{1} by 
\be
    \IE\{g_h(W+1) - W g_h(W)\},
\ee
where $g_h$ is the solution to the Stein equation
\ben                                                                \label{2}
	\lambda g_h(j+1) - j g_h(j) = h(j) - \IE h(Z), \quad j\in\IZ_+.
\ee
The left hand side of \eq{2} is an operator that characterises the Poisson
distribution; that is, for $\%A g(j) := \lambda g(j+1) - jg(j)$, 
\be
	\text{$\IE\%A g(Y) = 0$~~for all bounded $g$}
     \quad\iff\quad Y\sim\Po(\lambda).
\ee
Assume for simplicity that $W$ has the same support as $\Po(\lambda)$. With
\eq{2}, we can now write \eq{1} as
\ben                                                                \label{3}
	\dtv\bklr{\law(W),\Po(\lambda)} = \sup_{h\in\%H_{\tv}}
	\babs{\IE\%A g_h(W)}.
\ee
It turns out that \eq{3} is often easier to bound than \eq{1}.

\cite{Barbour1983} and \cite{Barbour1992} showed that, for all functions
$h\in\%H_\tv$,
\ben	                                                            \label{4}
	\norm{g_h} \leq 1\wedge\sqrt{\frac{2}{\lambda e}},
     \qquad \norm{\D g_h}\leq \frac{1-e^{-\lambda}}{\lambda},
\ee
where $\norm{\cdot}$ denotes the supremum norm and $\D g(j) := g(j+1)-g(j)$. So
here, if one is interested in the asymptotic $\lambda\toinf$, the Stein factors
are of order $\lambda^{-1/2}$ and $\lambda^{-1}$, respectively. With this we
have finished the first main step of Stein's method. 

As an example for the second step and also as a motivation for the main part of
this paper, assume that $W$ is a non-negative integer-valued random variable and
assume that $\tau$ is a function such that 
\ben                                                                \label{5}
    \IE\bklg{(W-\lambda)g(W)} = \IE\bklg{\tau(W)\D g(W)}
\ee
for all bounded functions $g$; see \cite{Cacoullos1994} and
\cite{Papathanasiou1995} for more details on this approach. To estimate the
distance between $\law(W)$ and the Poisson distribution with mean $\lambda$, we
simply use \eq{3} in connection with \eq{5} to obtain
\besn                                                                \label{6}
    \dtv\bklr{\law(W),\Po(\lambda)}
    &= \sup_{h\in\%H_\tv}\babs{\%A g_h(W)} \\
    &= \sup_{h\in\%H_\tv}\babs{\IE\bklg{\lambda g_h(W+1) - W g_h(W)}} \\
    &=\sup_{h\in\%H_\tv}\babs{\IE\bklg{\lambda\D g_h(W) - (W-\lambda)g_h(W)}}\\
    &= \sup_{h\in\%H_\tv}\babs{\IE\bklg{(\lambda-\tau(W))\D g_h(W)}} \\
    &\leq \frac{1-e^-\lambda}{\lambda}\IE\babs{\tau(W)-\lambda},
\ee
where for the last step we used \eq{4}. Thus, \eq{6} expresses the
$\dtv$-distance between $\law(W)$ and $\Po(\lambda)$ in terms of the average
fluctuation of $\tau$ around~$\lambda$. It is not difficult to show that $\tau
\equiv
\lambda$ if and only if $W \sim \Po(\lambda)$. 

Assume now that, for a fixed positive integer
$k$, $\tau(w) = \lambda + \delta_k(w)$, where $\delta_k(w)$ is the Kronecker
delta, and assume that $W_k$ is a random variable satisfying \eq{5} for
this~$\tau$. In this case we can in fact replace the last inequality in \eq{6}
by an equality to obtain
\ben                                                                \label{7}
    \dtv\bklr{\law(W_k),\Po(\lambda)} =
        \IP[W_k=k]\sup_{h\in\%H_\tv}\abs{\D g_h(k)}.
\ee
From Eq.~(1.22) of the proof of Lemma 1.1.1 of \cite{Barbour1992}
we see that, for $k=\floor{\lambda}$,
\ben                                                                \label{8}
    \sup_{h\in\%H_\tv}\abs{\D g_h(k)} \asymp \lambda^{-1}
\ee
as $\lambda\toinf$. Thus, \eq{7} gives
\ben                                                                \label{9}
    \dtv\bklr{\law(W_k),\Po(\lambda)} \asymp
        \IP[W_k=k]\lambda^{-1}
\ee
as $\lambda \toinf$. Note that, irrespective of the order of $\IP[W_k=k]$, the
asymptotic
\eq{9} makes full use of the second Stein factor of \eq{4}. To see that
$\law(W_k)$ in
fact exists, we rewrite \eq{5} as $\IE\%B_k g(W_k) = 0$, where
\besn                                                                \label{10}
    \%B_k g(w) & = \%A g(w) + \delta_k(w)\D g(w) \\
    & = \bklr{\lambda+\delta_k(w)}g(w+1) - \bklr{w+\delta_k(w)}g(w).
\ee
Recall from \cite{Barbour1988}, that $\%A$ can be
interpreted as the generator of a Markov process; in our case, as an
immigration-death process, with immigration rate $\lambda$, per capita death
rate~$1$  and $\Po(\lambda)$ as its stationary distribution. Likewise, we can
interpret $\%B_k$ as a perturbed immigration-death process with the same
transition rates, except in point~$k$, where the immigration rate is increased
to $\lambda+1$ and the per capita death rate is increased to $1+1/k$. Thus,
$\law(W_k)$ can be seen as the stationary distribution of this perturbed
process.

If $k=\floor{\lambda}$, the perturbation of the transition rates at point
$k$ is of smaller order than the transition rates of the corresponding
unperturbed immigration-death process in $k$. Thus, heuristically, 
$\IP[W_k=k]$ is of the same order as the probability $\Po(\lambda)\{k\}$ of the
stationary distribution of the unperturbed process, hence $\IP[W_k=k]
\asymp\lambda^{-1/2}$, and \eq{9} is of
order $\lambda^{-3/2}$. We omit a rigorous proof of this statement.

\begin{remark}
Note that by rearranging \eq{7} we obtain
\ben                                                    \label{14}
	\sup_{h\in\%H_{\tv}}\babs{\D g_h(k)}
     = \frac{\dtv\bklr{\law(W_k),\law(Z)}}{\IP[W_k = k]}.
\ee
for positive $k$. We can assume without loss of generality that
$g_h(0)=g_h(1)$ for all test functions $h$ because the value of $g_h(0)$ is not
determined by \eq{2} and can in fact be arbitrarily chosen. Thus $\D g_h(0) = 0$
and, taking the supremum over all $k\in\IZ_+$, we obtain
\ben                                                    \label{15}
    \sup_{h\in\%H_{\tv}}\norm{\D g_h}
     = \sup_{k\geq1}\frac{\dtv\bklr{\law(W_k),\law(Z)}}{\IP[W_k = k]}.
\ee
This provides a new interpretation of the bound $\norm{\D g_h}$ (a similar
statement can be made for $\norm{g_h}$, but then with a different family of
perturbations), namely as the ratio of the total variation distance between
some very specific perturbed Poisson distributions and the Poisson distribution,
and the probability mass at the location of these perturbations.

Let us quote \cite{Chen1998}, page 98:

\begin{quote}
\small Stein's method may be regarded as a method of constructing certain
kinds
of
identities which we call Stein identities, and making comparisons between
them. In applying the method to probability approximation we construct two
identities, one for the approximating distribution and the other for the
distribution to be approximated. The discrepancy between the two distributions
is then measured by comparing the two Stein identities through the use of the
solution of an equation, called Stein equation. To effect the comparison,
bounds on the solution and its smoothness are used. 
\end{quote}

Equations \eq{14} and \eq{15} make this statement precise. They express how
certain elementary deviations from the Stein identity of the approximating
distribution will influence the distance of the resulting distributions in the
specific metric, and they establish a simple link to the properties of the
solutions to \eq{2}. We can thus see $W$ from \eq{5} as a `mixture' of such
perturbations which is what is effectively expressed by Estimate~\eq{6}. 
\end{remark}

Thus, to understand why in some of the applications the Stein factors are not
as satisfying as in the above Poisson example, we will in the following sections
analyse the corresponding perturbed distributions in the cases of multivariate
Poisson and Poisson point processes.

In order to define the perturbations
to obtain an equation of the form \eq{7}, some care is needed, though. The
attempt to simply add the perturbation as in \eq{10}, may lead to an operator
that is not interpretable as the generator of a Markov process and thus the
existence of the perturbed distribution would not be guaranteed as easily.
It turns out that with suitable symmetry assumptions we can circumvent this
problem.

\section{Multivariate Poisson distribution} \label{sec2}

Let $d\geq 2$ be an integer, $\mu = (\mu_1,\dots,\mu_d)\in\IR_+^d$ such that
$\sum \mu_i=1$, and let $\lambda>0$.
Let $\Po(\lambda\mu)$ be the distribution on $\IZ_+^d$ defined as
$\Po(\lambda\mu) = \Po(\lambda\mu_1)\otimes\dots\otimes \Po(\lambda\mu_d)$.
Stein's method for multivariate Poisson approximation was introduced by
\cite{Barbour1988}; but 
see also \cite{Arratia1989}. Let $\eps^{(i)}$ denote $i$th unit vector. Using
the Stein operator 
\be
	\%A g (w) := \sum_{i=1}^d \lambda\mu_i\bklg{g(w+\eps^{(i)}) - g(w)}
     + \sum_{i=1}^d w_i\bklg{g(w-\eps^{(i)}) - g(w)}
\ee
for $w\in\IZ_+^d$, it is proved in
Lemma~3 of \cite {Barbour1988} that the solution $g_A$ to the Stein
equation $\%A g_A (w) = \delta_A(w)- \Po(\lambda\mu)\{A\}$ for
$A\subset\IZ_+^d$, satisfies
the bound
\ben										
			        \label{16}
	\bbbnorm{\sum_{i,j=1}^d \alpha_i\alpha_j \D_{i j} g_A}\leq
	\min\bbbklg{\frac{1 + 2\log^+(2\lambda)}{2\lambda}
    \sum_{i=1}^d\frac{\alpha_i^2}{\mu_i},\sum_{i=1}^d \alpha_i^2}   
\ee
for any $\alpha\in\IR^d$, where
\be
	\D_{i j}g(w) := g(w+\eps^{(i)} +\eps^{(j)})
                    -  g(w+\eps^{(i)}) - g(w+\eps^{(j)}) + g(w).
\ee

Let now $m_i = \floor{\lambda\mu_i}$ for $i=1,\dots,d$ and define 
\ben                                                            \label{18}
    A_1 = \{w\in \IZ_+^d : 0\leq w_1 \leq m_1, 0\leq w_2 \leq m_2\}.
\ee
\cite{Barbour2005} proves that, if $\mu_1,\mu_2>0$ and $\lambda\geq
(e/32\pi)(\mu_1\wedge\mu_2)^{-2}$, then
\ben										
				\label{19}
	\babs{\D_{12} g_{A_1}(w)} \geq
    \frac{\log\lambda}{20\lambda\sqrt{\mu_1\mu_2}}
\ee
for any $w$ with $(w_1,w_2) = (m_1,m_2)$. It is in fact not difficult to see
from the proof of \eq{19} that this bound also holds for the other quadrants
having corner~$(m_1,m_2)$.

\begin{figure}[t]
\center\includegraphics[width=0.8\hsize]{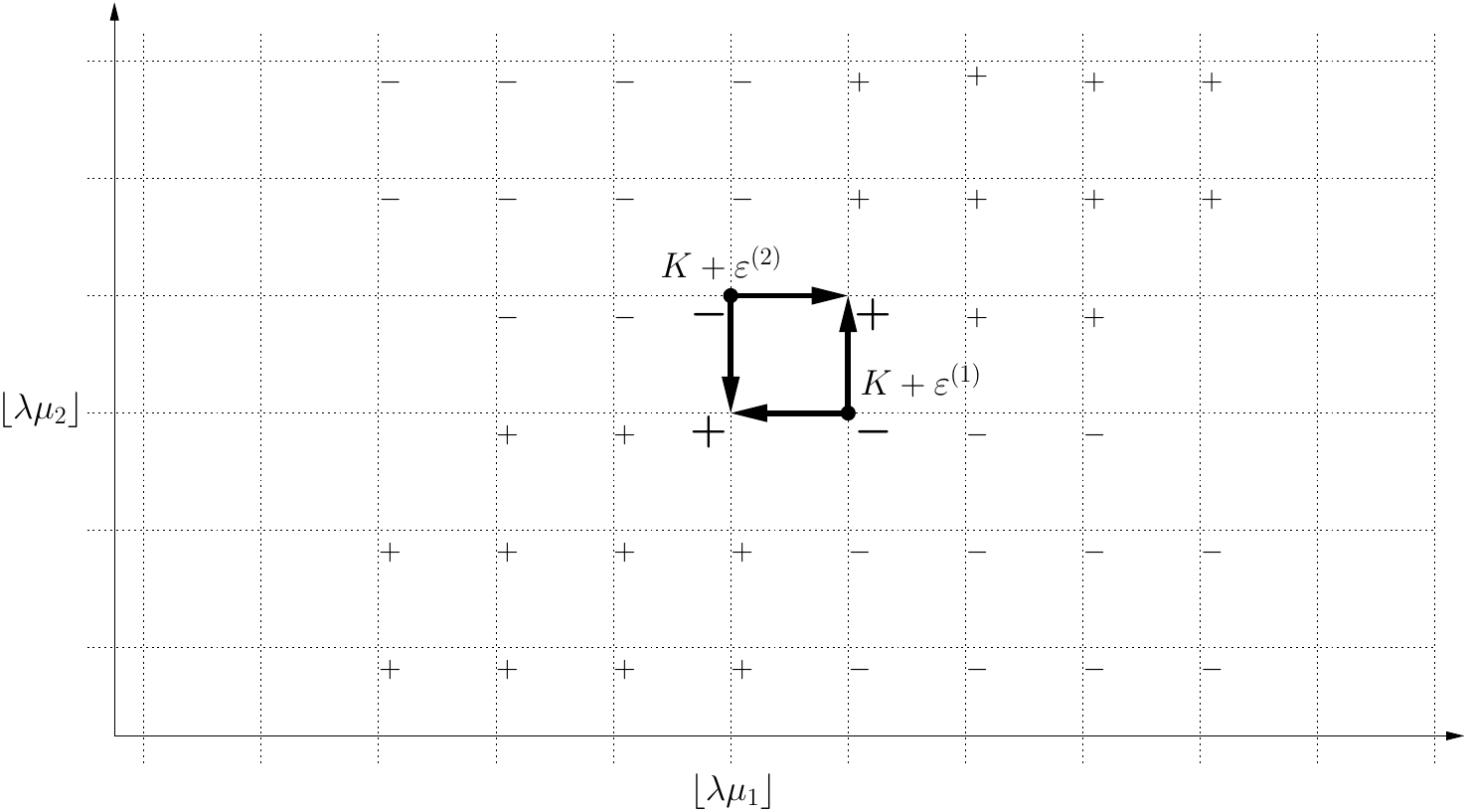}
\caption{\label{fig2} A rough illustration of the perturbed process defined by
the generator (\ref{20}). Between any of two connected points on the lattice
$\IZ_+^2$, we assume the transition dynamics of a unperturbed immigration-death
process, that is, in each coordinate immigration rate $\lambda\mu_i$ and per
capita death rate $1$. The arrows symbolise the additional perturbations with 
respect to the unperturbed immigration-death process; each arrow indicates an
increase by $1/2$ of the corresponding transition rate. The resulting
differences of the point probabilities between the equilibrium distributions of
the perturbed and unperturbed processes are indicated by the symbols $+$ and
$-$. The corresponding signs in each of the quadrants are heuristically
obvious, but they can be verified rigorously using the Stein equation,
Eq.~\eq{22}, and Eq.~$(2.8)$ of \cite{Barbour2005}.}
\end{figure}

\begin{example} \label{ex1}
Assume that $W$ is a random vector having the equilibrium distribution of the
$d$-dimensional birth-death process with generator  
\besn								\label{20}
	\%B_K g(w)  & =  \%Ag(w) \\
      &\quad +\ahalf\delta_{K+\eps^{(2)}}(w)\bkle{g(w+\eps^{(1)}) - g(w)}
      +\ahalf\delta_{K+\eps^{(1)}}(w)\bkle{g(w+\eps^{(2)}) - g(w)}\\
    &\quad +\ahalf\delta_{K+\eps^{(1)}}(w)\bkle{g(w-\eps^{(1)}) - g(w)} 
       +\ahalf\delta_{K+\eps^{(2)}}(w)\bkle{g(w-\eps^{(2)}) - g(w)} \\
     & = 
    \sum_{i=1}^d \bklr{\lambda\mu_i 
    + \ahalf\delta_1(i)\delta_{K+\eps^{(2)}}(w)
    +\ahalf\delta_2(i)\delta_{K+\eps^{(1)}}(w)}\bkle{g(w+\eps^{(i)}) - g(w)}\\
    &\quad + \sum_{i=1}^d
    \bklr{w_i+\ahalf\delta_1(i)\delta_{K+\eps^{(1)}}(w)
      +\ahalf\delta_2(i)\delta_{K+\eps^{(2)}}(w)}\bkle{g(w-\eps^{(i)}) - g(w)}
\ee
where $K=(m_1,m_2,\dots,m_d)$. Assume further that $\mu_1 = \mu_2$ , thus $m_1 =
m_2$ (the `symmetry condition'). See Figure~\ref{fig2} for an illustration of
this process. As the perturbations are symmetric in the first coordinates the
stationary distribution will also be symmetric in the first two
coordinates. 


Now, noting that for any bounded $g$ we have $\IE \%B_K g(W)=0$,
\besn                                                           \label{22}
    \IE\%A g(W)  
    & = \IE\%Ag(W) - \IE\%B_K g(W)\\
    & = - \ahalf\IP[W=K+\eps^{(2)}] \bkle{g(K+\eps^{(2)}+\eps^{(1)})-
	  g(K+\eps^{(2)})}\\
    &\quad- \ahalf\IP[W=K+\eps^{(1)}] \bkle{g(K+\eps^{(1)}+\eps^{(2)})-
	   g(K+\eps^{(1)})}\\    
    &\quad- \ahalf\IP[W=K+\eps^{(1)}]\bkle{g(K) - g(K+\eps^{(1)})} \\
    &\quad- \ahalf\IP[W=K+\eps^{(2)}]\bkle{g(K) - g(K+\eps^{(2)})} \\
    & = -   \IP[W=K+\eps^{(1)}]\D_{12}g(K),
\ee
where we used $\IP[W=K+\eps^{(1)}]=\IP[W=K+\eps^{(2)}]$ for the last
equality. Thus
\bes
	\dtv\bklr{\law(W),\Po(\lambda\mu)} & =
    \sup_{h\in\%H_\tv}\babs{\IE\%Ag_h(W)}\\
	& = \IP[W=K+\eps^{(1)}]\sup_{h\in\%H_\tv}\abs{\D_{12}g_h(K)}\\
	& \geq
\frac{\IP[W=K+\eps^{(1)}]\log\lambda}{20\lambda\sqrt{\mu_1\mu_2}}.
\ee
On the other hand, from \eq{16} for $\alpha=\eps^{(1)}$, $\alpha=\eps^{(2)}$ and
$\alpha=\eps^{(1)}+\eps^{(2)}$ respectively, it follows that 
\be
	\abs{\D_{12}g_h(w)} \leq
    \frac{\bklr{1+2\log^+(2\lambda)}(\mu_1+\mu_2)}{2\lambda\mu_1\mu_2}.
\ee
This yields the upper estimate
\bes
	\dtv\bklr{\law(W),\Po(\lambda\mu)} & =
\IP[W=K+\eps^{(1)}]\sup_{h\in\%H_\tv}\abs{\D_{12}g_h(K)} \\
	&\leq
\IP[W=K+\eps^{(1)}]\frac{\bklr{1+2\log^+(2\lambda)}(\mu_1+\mu_2)}{
2\lambda\mu_1\mu_2},
\ee
and thus we finally have
\ben                                                \label{23}
    \dtv\bklr{\law(W),\Po(\lambda\mu)}
    \asymp
    \frac{\IP[W=K+\eps^{(1)}]\log\lambda}{\lambda}.
\ee
Now, again heuristically, $\IP[W=K+\eps^{(1)}]$ will be of the order
$\Po(\lambda\mu_1)\klg{m_1}
\times\cdots\times\Po(\lambda\mu_d)\klg{m_d}\asymp\lambda^{d/2
} $ , so that \eq{23} will be of order $\log\lambda/\lambda^{1+d/2}$.

Recalling that the test function \eq{18} and also the corresponding test
functions for the other three quadrants are responsible for the logarithmic
term in \eq{23}, we may conclude a situation as illustrated
in Figure~\ref{fig2} for $d=2$. Different to the one-dimensional case,
where the perturbation moves probability mass from the point of the perturbation
to the rest of the support in a uniform way, the perturbations of the form
\eq{20} affect the rest of the support in a non-uniform way. However, further
analysis is needed to find the exact distribution of the probability
mass differences within each of the quadrants.

Note that the perturbation \eq{20} is `expectation neutral', that is, $W$ has
also expectation $\lambda\mu$, which can be seen by using $\IE\%B g(W)=0$ with
the function $g_i(w) = w_i $ for each coordinate $i$.

\end{example}

\section{Poisson point processes}

Stein's method for Poisson point process approximation
was derived by \cite{Barbour1988} and \cite{Barbour1992b}. They use the
Stein operator
\be
    \%A g(\xi) = 
    \int_{\Gamma}\bkle{g(\xi+\delta_\alpha) - g(\xi)} \lambda(d\alpha)
    +\int_{\Gamma}\bkle{g(\xi-\delta_\alpha) - g(\xi)} \xi(d\alpha),
\ee
where $\xi$ is a point configuration on a compact metric space $\Gamma$ and
$\lambda$ denotes the mean measure of the process. The most successful
approximation results have been
obtained in the so-called $d_2$-metric; see for example \cite{Barbour1992b},
\cite{Brown2000} and \cite{Schuhmacher2005}. Assume that $\Gamma$ is equipped
with a metric $d_0$ which is, for convenience, bounded
by $1$. Let $\%F$ be the set of functions $f:\Gamma\to\IR$,
satisfying
\be
    \sup_{x\neq y\in\Gamma} \frac{\abs{f(x)-f(y)}}{d_0(x,y)}\leq 1.
\ee
Define the metric $d_1$ on the set of finite measures on $\Gamma$ as
\be
    d_1(\xi,\eta) = 
    \begin{cases}
    1 & \text{if $\xi(\Gamma)\neq\eta(\Gamma)$,}\\
    \xi(\Gamma)^{-1}\sup\limits_{f\in\%F}\bbabs{\int f d\xi-\int f d\eta} 
    & \text{if $\xi(\Gamma)=\eta(\Gamma)$.}
    \end{cases}
\ee
Let now $\%H_2$ be the set of all functions from the set of finite measures into
$\IR$ satisfying
\be
    \sup_{\eta\neq\xi} \frac{\abs{h(\eta)-h(\xi)}}{d_1(\xi,\eta)}\leq 1.
\ee
We then define for two random measures $\Phi$ and $\Psi$ on $\Gamma$ the
$d_2$-metric as
\be
    d_2\bklr{\law(\Phi),\law(\Psi)}
    := \sup_{h\in\%H_2}\babs{\IE h(\Phi) - \IE h(\Psi)};
\ee 
for more details on the $d_2$-metric see \cite{Barbour1992} and
\cite{Schuhmacher2008}.

If $h\in\%H_2$ and $g_h$ solves the Stein
equation $\%Ag_h(\xi) = h(\xi) - \Po(\lambda)h$, \cite{Barbour1992b} prove the
uniform bound
\ben                                                        \label{24}
    \norm{\D_{\alpha\beta}g_h(\xi)} 
    \leq 1\wedge
\frac{5}{2\abs{\lambda}}\bbbklr{1+2\log^+\bbbklr{\frac{2\abs{\lambda}}{5}}},
\ee
where $\abs{\lambda}$ denotes the $L_1$-norm of $\lambda$ and where
\be
    \D_{\alpha\beta}g(\xi) =
        g(\xi+\delta_\alpha+\delta_\beta)-g(\xi+\delta_\beta)
        -g(\xi+\delta_\alpha)+g(\xi).
\ee
It has been shown by \cite{Brown1995} that the $\log$-term in \eq{24} is
unavoidable. However, \cite{Brown2000} have shown that it is possible to obtain
results without the~$\log$  using a non-uniform bound
on~$\D_{\alpha\beta}g_h$. 

Following the construction of \cite{Brown1995}, assume that $\Gamma
= S\cup\{a\}\cup\{b\}$,  where $S$ is a compact metric space,  $a$ and $b$ are
two additional points  with $d_0(a,b)= d_0(b,x) = d_0(a,x) = 1$
for all $x\in S$. Assume further that the measure $\lambda$ satisfies 
$\lambda(\{a\})=\lambda(\{b\}) = 1/\abs{\lambda}$ (again, the `symmetry
condition') and thus $\lambda(S) =
\abs{\lambda}-2/\abs{\lambda}$. For $m_a,m_b\in\{0,1\}$, define now  the test
functions 
\ben                                                        \label{25}
    h(\xi) = 
    \begin{cases}
     \frac{1}{\xi(\Gamma)} & \text{if $\xi(\{a\})=m_a$, $\xi(\{b\})=m_b$,
$\xi\neq 0$}\\
    0 & \text{else.}
    \end{cases}
\ee
It is shown by direct verification that $h\in\%H_2$. \cite{Brown1995} prove
that, for $m_a=m_b=1$, the corresponding solution $g_h$ to the Stein equation
satisfies the asymptotic
\ben                                                        \label{26}
    \abs{\D_{ab}g_h(0)} \asymp \frac{\log\abs{\lambda}}{\abs{\lambda}}  
\ee
as $\abs{\lambda}\toinf$, so that \eq{24} is indeed sharp, but it is easy
to see from their proof that \eq{26} will hold for the other
values of $m_a$ and $m_b$, as well.

\begin{figure}[t]
\center\includegraphics[width=0.5\hsize]{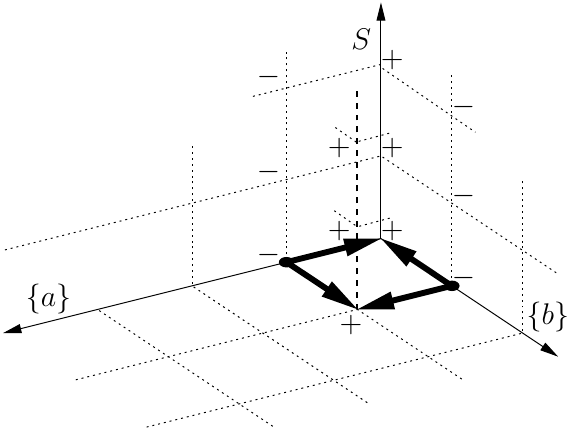}
\caption{\label{fig3} Illustration of the perturbed process
defined by the generator \eq{28} using the same conventions as in
Figure~\ref{fig2}. The corresponding signs can be obtained through the Stein
equation, Eq.~\eq{29} and the representation of the solution of the Stein
equation as in \cite{Brown1995}, for the different test functions \eq{25}.}
\end{figure}

\begin{example} \label{ex2}
Let $\Gamma$ and $\lambda$ be as above with the simplifying assumption that
$S$ is finite. Let $\Psi$ be a random point
measure with equilibrium
distribution of a Markov process with generator
\besn                                                           \label{28}
    \%B_0 g(\xi) 
    & = \%A g(\xi) + \ahalf\delta_{\delta_a}(\xi)\bkle{g(\xi+\delta_b)-g(\xi)}
      + \ahalf\delta_{\delta_b}(\xi)\bkle{g(\xi+\delta_a)-g(\xi)} \\
    &\quad + \ahalf\delta_{\delta_a}(\xi)\bkle{g(\xi-\delta_a)-g(\xi)}
	+ \ahalf\delta_{\delta_b}(\xi)\bkle{g(\xi-\delta_b)-g(\xi)}\\
    &= \int_{\Gamma}\bkle{g(\xi+\delta_\alpha)- g(\xi)}
    \bklr{\lambda+\ahalf\delta_{\delta_a}(\xi)\delta_b
	  +\ahalf\delta_{\delta_b}(\xi)\delta_a} (d\alpha)\\
    &\quad+\int_{\Gamma}\bkle{g(\xi-\delta_\alpha) - g(\xi)}
	\bklr{\xi+\ahalf\delta_{\delta_a}(\xi)\delta_a
	+\ahalf\delta_{\delta_b}(\xi)\delta_b}(d\alpha).
\ee
See Figure~\ref{fig3} for an illustration of this process.

Note that the situation here is different than in  Section~\ref{sec2}.
Firstly, we consider a weaker metric and, secondly, we impose a different
structure on~$\lambda$. Where as in Section~\ref{sec2} we assumed that the mean
of each coordinate is of the same order $\abs{\lambda}$, we assume now that
there are two special points $a$ and $b$ with $\lito(|\lambda|)$ mass attached
to them.
Again, in order to obtain a stationary distribution that is symmetric with
respect to $a$ and $b$, we impose the condition that the
immigration rates at the two coordinates $a$ and $b$ are the same.

Now, for any bounded function $g$,
\besn                                                               \label{29}
    \IE\%A g(\Psi)  & = \IE\%Ag(\Psi) - \IE\%B_0 g(\Psi)\\
    & = - \ahalf\IP[\Psi=\delta_b] 
            \bkle{g(\delta_a + \delta_b) - g(\delta_a)}
	- \ahalf\IP[\Psi=\delta_a] 
            \bkle{g(\delta_b + \delta_a) - g(\delta_b)}\\
    &\quad - \ahalf\IP[\Psi=\delta_a] 
            \bkle{g(\delta_a) - g(0)}
	- \ahalf\IP[\Psi=\delta_b] 
            \bkle{g(\delta_b) - g(0)}\\ 
    & = - \IP[\Psi = \delta_a]\D_{ab}g(0),
\ee
where we used that $\IP[\Psi = \delta_a]=\IP[\Psi = \delta_b]$.
Thus, using \eq{26},
\ben                                                               \label{30}
    d_2\bklr{\law(\Psi),\Po(\lambda)} 
    = \IP[\Psi=\delta_a]\sup_{h\in\%H_2}\abs{\D_{ab} g_h(0)} \asymp
    \frac{\IP[\Psi=\delta_a]\log\abs{\lambda}}{\abs{\lambda}}.
\ee

Figure~\ref{fig3} illustrates the situation for $\abs{\Gamma}=3$. If the process
$\Phi_t$ is somewhere on the bottom plane, that is $\Phi(S) = 0$, it will most
of the times quickly jump upwards, parallel to the $S$-axis, before
jumping between the parallels, as the immigration rate into $S$ is far
larger than the jump rates between the parallels. Thus, because of the
perturbations, probability mass is moved---as illustrated in
Figure~\ref{fig3}---not only between the
perturbed points but also between the parallels. Although indicator functions
are not in $\%H_2$, the test functions in \eq{25} decay slowly enough to
detect this difference.

\end{example}

\begin{remark}
Note again, as in Example~\ref{ex1}, that the perturbation in the above example
is neutral with respect to the measure $\lambda$. It is also interesting to
compare the total number of points to a Poisson distribution with mean
$\abs{\lambda}$ in the $\dtv$-distance. Note that \eq{29} holds in particular
for functions $g_h$ which depend only on the number of points of $\Psi$. Thus,
using \eq{3} in combination with \eq{29} yields
\be
    \dtv\bklr{\law(\abs{\Psi}),\Po(\abs{\lambda})} 
    = \IP[\Psi = \delta_a]\sup_{h\in\%H_\tv}\abs{\D^2 g_h(0)}
     \asymp \frac{\IP[\Psi = \delta_a]}{\abs{\lambda}},
\ee
where $\D^2 g(w) = \D g(w+1) - \D g(w)$ (which corresponds to
the first difference in~\eq{4}) and where we used the fact that $\abs{\D^2
g_h(0)}\asymp\abs{\lambda}^{-1}$, again obtained from the proof of
Lemma~1.1.1 of \cite{Barbour1992}. Thus we have effectively constructed an
example, where the attempt to match not only the number but also the location
of the points introduces an additional factor $\log\abs{\lambda}$ if using the
$d_2$-metric. 

\end{remark}

\section*{Acknowledgements}

The author would like to thank Gesine Reinert and Dominic Schuhmacher for
fruitful discussions and the anonymous referee for helpful comments.

\def\url#1{#1}


\begin{thebibliography}{19}
\providecommand{\natexlab}[1]{#1}

\bibitem{Arratia1989}
Arratia, R., Goldstein, L., Gordon, L.: 
Two moments suffice for {P}oisson approximations: the {C}hen-{S}tein
method. 
Ann. Probab. \textbf{17}, 9--25  (1989)

\bibitem{Barbour1988}
Barbour, A.D.: 
Stein's method and {P}oisson process convergence.
J. Appl. Probab. \textbf{25A}, 175--184  (1988)

\bibitem{Barbour2005}
Barbour, A.D.:
Multivariate {P}oisson-binomial approximation using {S}tein's method.
In: Stein's method and applications, Vol.~5 of Lect.
  Notes Ser. Inst. Math. Sci. Natl. Univ. Singap., pp.~131--142. Singapore
  University Press, Singapore (2005)

\bibitem{Barbour1992b}
Barbour, A.D., Brown, T.C.:
Stein's method and point process approximation.
Stochastic Process. Appl. \textbf{43}, 9--31  (1992)

\bibitem{Barbour1983}
Barbour, A.D., Eagleson, G.K.:
Poisson approximation for some statistics based on exchangeable
trials.
Adv. in Appl. Probab. \textbf{15}, 585--600  (1983)

\bibitem{Barbour1998a}
Barbour, A.D., Utev, S.A.:
Solving the {S}tein equation in compound {P}oisson approximation.
Adv. in Appl. Probab. \textbf{30}, 449--475  (1998)

\bibitem{Barbour1992c}
Barbour, A.D., Chen, L.H.Y., Loh, W.-L.:
Compound {P}oisson approximation for nonnegative random variables via
  {S}tein's method.
Ann. Prob. \textbf{20},  1843--1866  (1992)

\bibitem{Barbour1992}
Barbour, A.D., Holst, L., Janson, S.:
Poisson approximation. 
Oxford University Press, New York (1992)

\bibitem{Brown1995}
Brown, T.C., Xia, A.:
On {S}tein-{C}hen factors for {P}oisson approximation.
Statist. Probab. Lett. \textbf{23}, 327--332 (1995)

\bibitem{Brown2000}
Brown, T.C., Weinberg, G.V., Xia, A.:
Removing logarithms from {P}oisson process error bounds.
Stochastic Process. Appl. \textbf{87}, 149--165. (2000).

\bibitem{Cacoullos1994}
Cacoullos, T., Papathanasiou, V., Utev, S.A.:
Variational inequalities with examples and an application to the
  central limit theorem.
Ann. Probab. \textbf{22}, 1607--1618 (1994)

\bibitem{Chen1998}
Chen, L.H.Y.:
Stein's method: some perspectives with applications.
In: Probability towards 2000, volume 128 of Lecture Notes in Statistics,
pp.~97--122. Springer, New York (1998)

\bibitem{Chen2009}
Chen, L.H.Y., R\"ollin, A.: 
Stein couplings for normal approximation. Preprint (2010). Available via arXiv,
\url{http://arxiv.org/abs/1003.6039} 

\bibitem{Chen2004a}
Chen, L.H.Y., Shao, Q.-M.:
Normal approximation under local dependence.
Ann. Probab. \textbf{32}, 1985--2028 (2004)

\bibitem{Papathanasiou1995}
Papathanasiou, V., Utev, S.A.:
Integro-differential inequalities and the {P}oisson approximation.
Siberian Adv. Math. \textbf{5}, 120--132 (1995)

\bibitem{Schuhmacher2005}
Schuhmacher, D.:
Upper bounds for spatial point process approximations.
Ann. Appl. Probab. \textbf{15}, 615--651 (2005)

\bibitem{Schuhmacher2008}
Schuhmacher, D., Xia, A.:
A new metric between distributions of point processes.
Adv. Appl. Probab. \textbf{40}, 651--672 (2008)

\bibitem{Stein1972}
Stein, C.:
A bound for the error in the normal approximation to the distribution
  of a sum of dependent random variables.
In: Proceedings of the Sixth Berkeley Symposium on Mathematical
  Statistics and Probability Vol.~II, pp.~583--602.
University California Press, Berkeley (1972)

\bibitem{Stein1986}
Stein, C.:
Approximate computation of expectations. Institute of Mathematical Statistics,
Hayward (1986)

\end{thebibliography}
\end{document}